\documentclass[12pt,oneside,a4paper]{article}
\usepackage[cp1251]{inputenc}
\usepackage{amsfonts,amssymb,latexsym,amsmath}
\usepackage[T2A]{fontenc}
\usepackage[english]{babel}
\allowdisplaybreaks[1]
\begin{document}
\setcounter{page}{1}
\pagestyle{plain}

\large
\vskip 2cm

\bigskip
\noindent UDC 517.54
\ \newline D.L.~STUPIN
\ \newline Tver state university
\smallskip
\ \newline {\bf THE SHARP ESTIMATES OF ALL INITIAL}
\ \newline {\bf TAYLOR COEFFICIENTS IN THE}
\ \newline {\bf KRZYZ'S PROBLEM}
\footnote[1] {\copyright\;\; D.L.~Stupin, 2010.}
\smallskip


{\sl For each $t>0,$ up to the number $n=N(t),$ the exact estimations of all
initial Taylor coefficients in the class $B_t$ were found, where $B_t$ is
a set of holomorphic in unit disk functions $f,$ $0<|f|<1,$ $f(0)=e^{-t}.$}

\smallskip
{\bf 1.}
Let $\Delta:=\{z\in\mathbb C\,:\,|z|<1\},$ where $\mathbb C$ --- a field of all
complex numbers.

Class $B$ consists of holomorphic in $\Delta$ functions $w=f(z),$ not vanishing
and such that $|f(z)|\leqslant1,$ $z\in\Delta.$

The Krzyz conjecture~[1] consists in that for every $f\in B$ and all
natural $n$ the Taylor coefficients $\{f\}_n$ of function $f$ satisfy to inequality
$|\{f\}_n|\leqslant2/e.$
The extremals can be only rotations of the function
$F^*(z^n,1)$ in the planes $z$ and $w,$ where
$$
  F^*(z,t):=e^{-t\,H(z)}, \qquad H(z):=\frac{1+z}{1-z}, \quad t>0.              \eqno(1)
$$

The existence of extremals of this problem is obvious, since after addition
of function $f(z)\equiv0$ to the class $B$ it becomes a family of functions,
compact in the topology of the locall uniform convergence.
However, extremal prob\-lems on the class $B$ are very complicated and,
at present, the hypothesis remains unconfirmed.

Since the class $B$ is invariant under rotations in the planes of the variables
$z$ and $w,$ than we can restrict ourselves to studying only of functions,
for which $0<\{f\}_0<1.$
Further, we can fix the parameter $t\in(0,+\infty)$ and set
$\{f\}_0=e^{-t};$
these subclasses we denote by $B_t.$

From geometrical considerations it is clear that every function of class
$B_t$ can be represented in the form
$$
  f(z)=F^*(\omega(z),t), \qquad \omega\in \Omega,                               \eqno(2)
$$
where class $\Omega$ consists of holomorphic in unit disk $\Delta$
functions $\omega(z)$ such, that $\omega(0)=0$ and $|\omega(z)|<1,$
when $z\in\Delta.$

Note that for every $t>0 $ the formula \thetag{2} establishes a one-to-one
cor\-res\-pon\-den\-ce between the classes $B_t$ and $\Omega$ (see [2]).

If $t\in(0,2]$ than the Krzyz conjecture can be refined:
if $f\in B_t$ than, probably, the accurate estimates
$|\{f\}_n|\leqslant|\{F^*\}_1(t)|={2t}/{e^t}$ are correct;
the equality is delivered only by functions
$F^*(e^{i\varphi}z^n,t),$ $n\in\mathbb N,$ $\varphi\in\mathbb R.$

\smallskip
{\bf 2.}
Let us dwell on representations of the form \thetag{2}.
And let the functions $g(z)$ and $G(z)$ be holomorphic in $\Delta.$
Function $g(z)$ is called subordinate in the disk $\Delta$ for the function $G(z),$
if it can be represented in $\Delta$ in the form
$g(z)=G(\omega(z)),$ where $\omega\in\Omega.$
Function $G(z)$ will be called majorant for the $g(z)$ in the domain $\Delta.$

Let
$g(z)=\sum\limits_{n=0}^\infty{{\{g\}_n}z^n},$
$G(z)=\sum\limits_{j=0}^\infty{{\{G\}_j}z^j},$
$\omega(z)=\sum\limits_{k=1}^\infty{{\{\omega\}_k}z^k}.$
Then
$$
  g(z)=G(\omega(z))=
  \sum\limits_{j=0}^\infty{\{G\}_j\omega(z)^j}=
  \{G\}_0+
  \sum\limits_{n=1}^\infty
  {
    \left(
    \sum\limits_{j=1}^n
    {
      \{G\}_j\{\omega^j\}_n
    }
    \right)
    z^n
  }.
$$
Whence it appears $\{g\}_0=\{G\}_0$ and
$$
  \{g\}_n(\omega)=\sum\limits_{j=1}^n{\{G\}_j\{\omega^j\}_n},
  \qquad n\in\mathbb N,
  \quad\omega\in\Omega.                                                         \eqno(3)
$$

\smallskip
{\bf 3.}
The class of all functions $f(z),$ regular and univalent in $\Delta,$ with
the normalization $f(0)=0,$ $f'(0)=1,$ mapping the disk $\Delta$ on convex
domain, is denoted by $S^0.$

The set of all functions $h(z),$ with a positive in $\Delta$ real part and
with the normalization $h(0)=1,$ mapping the unit disk $\Delta$ into
right half-plane, is called the Caratheodory class and is denoted by $C.$

Between $C$ and $S^0$ there is the bijection [2]:
$$
  h(z)=1+z\frac{f''(z)}{f'(z)}, \qquad h\in C, \quad f\in S^0.                  \eqno(4)
$$

The following simple but important statement is hold~[2]:
\newtheorem{Lem}{Lemma}
\begin{Lem}
If the function $s(z)=\sum\limits_{n=1}^\infty{\{s\}_nz^n},$ regular in
$\Delta,$ is subordinated to the function $S(z)\in S^0,$
then the sharp estimates
$$
  |\{s\}_n|\leqslant|\{S\}_1|=1, \qquad n\in\mathbb N,
$$
are valid.
Equality is achieved only on the functions
$S(e^{i\varphi}z^n),$ $\varphi\in\mathbb R,$ $n\in\mathbb N.$
\end{Lem}

Carath\'eodory and T\"oplitz [3,~4] have completely solved the problem about
the possibility of extention of a polynomial up to a function of the class $C.$
Later Shur gave a constructive proof of this theorem [5].

\newtheorem{Criterion}{Criterion}
\begin{Criterion}[Carath\'eodory,\,T\"oplitz]
Let $\{h\}_1,\ldots,\{h\}_n,$ $n\in\mathbb N,$ --- are fixed complex numbers.
Polynom
$$
  R(z,n):=1+\sum\limits_{k=1}^n{\{h\}_kz^k}
$$
can be extended up to the function
$h(z)= R(z,n)+O(z^{n+1})\in C$ if and only if determinants
$$
  M_k:=\det{\{a_{ij}\}_{i,j=0}^{k}}, \qquad 1\leqslant k\leqslant n,
$$
$$
  a_{ii}=2, \qquad
  a_{ij}=\{h\}_{j-i}, \quad j>i, \qquad
  a_{ij}=\overline{a_{ji}}, \quad j<i,
$$
either all positive or positive to a certain number,
from which are equal to zero. In the latter case, the extension is unique.
\end{Criterion}

\smallskip
{\bf 4.}
For further progress we need to study the Taylor coefficients of our majorant function
$F^*(z, t)$ of class $B_t.$
Everywhere below, we will not be interested in the zero coefficient of this
function, since it is not included in the formula \thetag{3}.
The first coefficient $\{F^*\}_1(t)$ of function $F^*(z, t)$ is equal to
$-2t/e^t.$
We normalize the function $F^*(z, t)$ so, that the first coefficient in its Taylor
expansion becomes equal to $1$. Let us introduce the notation
$$
  F(z, t):= \dfrac{F^*(z, t)}{\{F^*\}_1(t)}.                                    \eqno(5)
$$

The question arises: whether convex univalent functions $f\in S^0$ exist there,
with some initial Taylor coefficients that match with all first coefficients
of the function $F(z, t),$ except for the coefficient $\{F\}_0(t)?$

It is not hard to check.
By substituting $f(z)=F(z,t)$ to the formula \thetag{4}, we obtain
$$
  h(z)=1+2z\left(\frac1{1-z}-\frac1{(1-z)^2}t\right).
$$
From which we elementary derive remarkably simple formula
$$
  \{h\}_j=2(1-jt), \qquad j\in\mathbb N.                                        \eqno(6)
$$

We use the Carath\'eodory-T\"oeplitz extension criterion of polynomials up to a
function of class $C.$
Let us compute the principal minors $M_{j-1},$ for all $j\in\mathbb N.$
Here the index $j-1$ means, that the dimension of the corresponding to the
minor $M_{j-1}$ matrix is equal to $j.$
According to lemma~2, the formulation and the proof of which can be found
at the end of this work (see section~6), we have:
$$
  M_{j-1}=2^{2(j-1)}t^{j-1}(2-(j-1)t), \qquad j\in\mathbb N.                    \eqno(7)
$$

Further, it is obvious that the minors $M_1,\ldots,M_{n-1}$ are not negative
if and only if $t\leqslant2/(n-1),$ for $n>1,$ and
$t>2,$ for $n=1,$ or $n\leqslant2/t+1.$
Note also, that if $t=2/(n-1),$ $n\in\mathbb N\setminus\{1\},$
the extension is unique.

Thus, we have

\newtheorem{Th}{Theorem}
\begin{Th}
For every $t>0$ and $n\leqslant2/t+1,$ $n\in\mathbb N,$ the segment of Taylor
expansion of the function $F(z,t),$ first introduced in the formula \thetag{5},
--- polynom $P(z, t, n):=z+\sum\limits_{k=2}^n{\{F\}_k(t)\,z^k}$ ---
can be extended to the function $f(z)= P(z,t,n)+O(z^{n+1})\in S^0.$
For $t=2/(n-1),$ $n\in\mathbb N\setminus\{1\},$ the extension is unique.
\end{Th}

\newenvironment{Proof}                 
{\par\noindent{\bf Proof.}}            
{}

\smallskip
{\bf 5.}
From lemma~1, theorem~1, formula~\thetag{3} and normalization~\thetag{5}
follows a central for this work

\begin{Th}
For every $t>0,$ arbitrary $N\leqslant2/t+1,$ $N\in\mathbb N,$
and each $f^*\in B_t,$ sharp estimates
$$
  |\{f^*\}_n|\leqslant|\{F^*\}_1(t)|=\dfrac{2t}{e^t},
  \qquad n\in\{1,\ldots,N\},                                                    \eqno(8)
$$
are correct.
Extremals in the estimates \thetag{8} are only the functions
$$
  F^*(e^{i\varphi}z^n,t), \qquad \varphi\in\mathbb R,
$$
where the function $F^*$ defined by \thetag{1}.
\end{Th}

\begin{Proof}
We fix $\omega\in\Omega,$ $t>0$ and $N\leqslant2/t+1,$ $N\in\mathbb N.$

Let us take a natural number $n,$ not exceeding the number $N.$
Using formula \thetag{3}, we write $n$-th coefficient of function
$$
  f(z):=F(\omega(z),t),
$$
where $F$ is defined in formula \thetag{5}, in the form
$$
  \{f\}_n=\sum\limits_{j=1}^n{\{F\}_j\{\omega^j\}_n}.
$$

Now we apply theorem~1 to $n$-th segment of Taylor expansion of function
$F(z,t),$ which we have denoted by $P(z,t,n).$
Let $S(z)$ --- be an extention of polynom $P(z,t,n)$ to function of class $S^0.$
Then, using the formula \thetag{3}, $n$-th coefficient of function
$$
  s(z):=S(\omega(z),t)
$$
can be written as
$$
  \{s\}_n=\sum\limits_{j=1}^n{\{S\}_j\{\omega^j\}_n}.
$$
From which, by lemma~1, we find that
$$
  |\{s\}_n|\leqslant1.
$$
But $\{S\}_j:=\{F\}_j(t),$ where $j\in\{1,\ldots,n\},$ therefore
$\{f\}_n=\{s\}_n,$ on basis of which we conclude, that
$$
  |\{f\}_n|\leqslant1.
$$
Remembering about the normalization \thetag{5}, we obtain
the estimates \thetag{8}.
Accuracy of the estimates \thetag{8} and the form of extremal functions
follow from lemma~1.

The theorem has been completely proved.
\end{Proof}

From theorem~2 it implies, that the smaller the number $t>0$ we fix, the more
Taylor coefficients we can estimate on the class $B_t.$
In this case, our estimates are sharp in the sense, that the equality,
in the inequality~\thetag{8}, is attained on functions
$F^*(e^{i\varphi}z^n,t).$

For example, if $t>2$ we can estimate only one coefficient, for $t=2$
--- two coefficients, for $t=1$ --- three coefficients, and at $t=1/2$
--- five. And so on.
Similar results were obtained in [6].

\smallskip
{\bf 6.}
We eliminate the blank in the arguments given above.
To do this we must only prove the validity of the formula \thetag{7}.

\begin{Lem}
If the coefficients $\{h\}_j,$ $j\in\mathbb N,$ are defined by formula \thetag{6},
then for all integers $n\geqslant0$
$$
  M_n=2^{2n}t^n(2-nt).
$$
\end{Lem}

\begin{Proof}
  The minor $M_n/2^{n+1}$ is equal to the determinant
  $$
  \begin{vmatrix}
  1          & 1-t      & 1-2t     &\ldots & 1-(n-1)t & 1-nt     \\
  1-t        & 1        & 1-t      &\ldots & 1-(n-2)t & 1-(n-1)t \\
  1-2t       & 1-t      & 1        &\ldots & 1-(n-3)t & 1-(n-2)t \\
  \vdots     & \vdots   & \vdots   &\ddots & \vdots   & \vdots   \\
  1-(n-1)t   & 1-(n-2)t & 1-(n-3)t &\ldots & 1-2t     & 1-t      \\
  1-nt       & 1-(n-1)t & 1-(n-2)t &\ldots & 1-t      & 1        \\
  \end{vmatrix}
  .
  $$
  Subtracting from each row, except the first one, the previous one we obtain
  $$
  \begin{vmatrix}
    1      & 1-t    & 1-2t   &\ldots & 1-(n-1)t & 1-nt   \\
    -t     & t      & t      &\ldots & t        & t      \\
    -t     & -t     & t      &\ldots & t        & t      \\
    \vdots & \vdots & \vdots &\ddots & \vdots   & \vdots \\
    -t     & -t     & -t     &\ldots & t        & t      \\
    -t     & -t     & -t     &\ldots & -t       & t      \\
  \end{vmatrix}
  .
  $$
  For each column, except for the last one, we add the last column
  $$
  \begin{vmatrix}
    2-nt   & 1-(n+1)t & 1-(n+2)t &\ldots & 1-(2n-1)t & 1-nt   \\
    0      & 2t       & 2t       &\ldots & 2t        & t      \\
    0      & 0        & 2t       &\ldots & 2t        & t      \\
    \vdots & \vdots   & \vdots   &\ddots & \vdots    & \vdots \\
    0      & 0        & 0        &\ldots & 2t        & t      \\
    0      & 0        & 0        &\ldots & 0         & t      \\
  \end{vmatrix}
  =
  $$
  $$
  =2^{n-1}t^n(2-nt).
  $$
\end{Proof}

\smallskip
{\bf 7.}
We present an example of an extention.
Let $t=1/2.$
By theorem~1, the desired extension is unique.
Setting in the formula \thetag{4}
$$
  f(z)=z+\frac12z^2+\frac16z^3-\frac1{24}z^4-\frac{19}{120}z^5+\ldots \in S^0
$$
or using the formula \thetag{6}, we obtain that
$$
  h(z)=1+z-z^3-2z^4+\ldots \in C.
$$
We know that the function
$$
  \omega(z)
  =\frac{1-h(z)}{1+h(z)}
  =-\frac12z+\frac14z^2+\frac38z^3+\frac9{16}z^4+\ldots
$$
belongs to the class $\Omega.$
It is also known (see~[2, 5]) that
$$
  \omega(z)=
  \lambda\frac
  {\overline\alpha_{n-1}+\overline\alpha_{n-2}z+\ldots+\overline\alpha_0z^{n-1}}
  {\alpha_0+\alpha_1z+\ldots+\alpha_{n-1}z^{n-1}}.
$$
In this case $\lambda=1$ (see~ [2, 5]).
Having found the parameters $\alpha_0,\ldots,\alpha_{n-1},$ we find
$$
  \omega(z)=z\frac{1-z^2-2z^3}{-2-z+z^3}.
$$
Whence
$$
  h(z)=\frac{1+z-z^3-z^4}{1+z^4}.
$$
Now we use the formula \thetag{4}, by substituting the
obtained expression for $h(z)$ there. Well
$$
  f(z)=\int\limits_0^z
  {
    \frac{\left(\frac{v^2+\sqrt2v+1}{v^2-\sqrt2v+1}\right)^{\sqrt2/4}}{\sqrt{1+v^4}}dv
  }
  .
$$

\bigskip
CITATION

\begin{enumerate}
  \item Krzyz~J.G.
  {\it Coefficient problem for bounded nonvanishing fun\-cti\-ons}~//
  Ann. Polon. Math. 1968. V.~70. P.~314.
  \item Golusin~G.M.
  {\it Geometric theory of functions of a complex variable.}
  Engl.~transl.: American Mathematical Society, Providence, RI,~1969.
  \item Carath\'eodory~C.
  {\it \"Uber die Variabilit\"atsbereich
  des Fourierschen Kon\-stan\-ten von Positiv Harmonischen Funktion}~//
  Rendiconti Circ. Mat. di~Palermo. 1911. V.~32. P.~193--217.
  \item T\"oplitz~O.
  {\it \"Uber die Fouriersche Entwicklung Positiver Funktionen}~//
  Rendiconti. Circ. Mat. di~Palermo. 1911. V.~32. P.~191--192.
  \item Schur~I.
  {\it \"Uber potenzreihen, die in Innern des Einheitskrises Beschr\"ankt Sind}~//
  Reine Angew. Math. 1917. V.~147. P.~205--232.
  \item Peretz~R.
  {\it Applications of subordination theory to the class of bounded nonvanishing
  functions}~//
  Compl. Var. 1992. V.~17. P.~213--222.
\end{enumerate}

\begin{flushright}
April 2010
\end{flushright}

\end{document}